\documentclass [11pt, oneside]{amsart}
\usepackage {amsfonts}
\usepackage {amsmath}
\usepackage {amsthm}
\usepackage {amssymb}
\usepackage{framed}
\usepackage {amsxtra}
\usepackage {enumerate}
\input xy
\xyoption{all}

\ifx\pdftexversion\undefined
  \usepackage[dvips]{graphics}
\else
  \usepackage[pdftex]{graphics}
\fi

\theoremstyle{definition}
\newtheorem{df}{Definition} [section]
\newtheorem{dfs}[df]{Definitions}
\newtheorem{ex}[df]{Example}
\newtheorem{exs}[df]{Examples}

\newtheorem{remark}[df]{Remark}
\newtheorem{remarks}[df]{Remarks}

\theoremstyle{plain}
\newtheorem{thm}[df]{Theorem}

\newtheorem{lemma}[df]{Lemma}

\newtheorem{prop}[df]{Proposition}

\title{A Calculus of Inconsistency I: Sentential Logic}
\author{Dan Seabold \and Stefan Waner \and Steve Warner}
\address{Dan Seabold: Department of Mathematics\\
       Hofstra University}
\email{Daniel.E.Seabold@Hofstra.edu}
\address{Stefan Waner: Department of Mathematics\\
       Hofstra University}
\email{Stefan.Waner@Hofstra.edu}
\address{Steve Warner: Department of Mathematics\\
       Hofstra University}
\email{Steve.Warner@Hofstra.edu}
\date {September 2004}
\begin{document}
\maketitle
\begin{abstract}
We describe a graph-theoretic syntax for self-referential formulas
as well as a four-valued logic to include contradictory and
independent formulas. We then explore the degree to which
generalized truth tables can be realized in our theory, and go on
to describe a model theory for sentential calculus, wherein models
are allowed to include contradictions (such as the ``Liar'') and
formulas that result from them as an integral part of their
structure. This sets the groundwork for a sequel in which we
construct models of set theory that include contradictions.
\end{abstract}

\newcommand{\dom}{{\rm dom \ }}
\newcommand{\ran}{{\rm ran \ }}

\section{Introduction}\label{Introduction}
Any logic that permits self-reference opens itself to paradox, as
illustrated by the Liar's paradox, a simple form of which states
``This statement is false.'' In their study of the subject,
Barwise and Etchemendy \cite{Barwise} address this and similar
paradoxes by introducing a ``non-Russellian'' logic that includes
a hierarchy of discourse in which self-referential statements,
such as the Liar, are either true or false.

If we wish to include self-referential statements without
abandoning the Russellian system of logic, we are faced with
contradictions. Classically, one has the high school proof that
any contradiction implies that all statements are contradictions,
thereby contaminating the whole universe with contradiction.
Consider the simplest variant of this argument: ``If $p$ is true
then $p\vee q$ is true for any formula $q$, but then if $p$ is
also false, it follows that $q$ must be true.'' In our analysis of
this argument, if $p$ is a contradiction and $q$ is false, then
$p\vee q$ is also a contradiction, and our rules for computing
truth will not permit us to conclude anything more about $q$. In
this way we avoid contaminating the entire model.

The approach we shall use to compute truth is the traditional
hierarchical approach: the truth of formulas that are not
subformulas of axioms can only be computed from ``below''; that
is, from a knowledge of the truth of subformulas. They can only be
inferred from above if they are subformulas of axioms, or through
arguments (such as rules of inference) outside the model.

With regard to rules of inference, we shall see that all the
classical rules of inference still apply, but that few of them
 are tautologies (in particular, modus ponens is not a tautology).

Our models of logic will be based on graphs analogous to those
used by Axcel to model set theory with antifoundation in
\cite{Axcel}. In this paper we restrict our syntax to sentential
logic; richer forms of syntax such as that in Smullyan
\cite{Smullyan} will be developed in a future paper. We shall also
find that the most natural setting in which to include
contradictions is through the use of a four-valued system of
logic, where the two additional values are $L$ (``is a lie") for
contradictions, and $V$ (``is vacuous") for independent
statements. In an unpublished manuscript, Linton \cite{Linton}
outlines a four-valued system of logic whose four truth values are
somewhat reminiscent of ours, although he does not interpret them
as we do here, nor does he develop a calculus of their use or an
application to self-referential statements.

In \S 2, we describe the basic calculus of our four-valued system
of logic and observe that all formulas in classical propositional
logic can be identified with special formulas in our sense. In \S
3, we describe (Proposition \ref{realization1} and Theorem
\ref{realization}) exactly what kinds of truth tables can be
realized by formulas. In \S 4 we outline our model theory for
sentential calculus, describe how to determine truth in a model,
and give models of sentential calculus that extend the classical
one; one in which the Liar is a contradiction, and another in
which the Liar is false.

The authors are grateful to Sylvia Silberger and Aileen Michaels
for many stimulating conversations, and especially to Scott Davis
for his Honors thesis work on computing all possible unary gates
and for his observations on an early version of this work.

\section{Preliminaries}\label{Preliminaries} Axcel \cite{Axcel}
constructed anti-foundational models of set theory using pointed
graphs with cycles to model membership hierarchies in sets. In the
same vein, the formulas in our version of sentential logic will be
pointed graphs with cycles. If $\Gamma$ is a set, recall that a
$\Gamma$-\textbf{graph} consists of a finite set $Z$ of {\bf
directed edges} together with functions $S: Z\to \Gamma$ and
$T: Z \to \Gamma$ called {\bf source} and {\bf target} maps
respectively. We refer to the sources and targets of the arrows as
\textbf{nodes.}  If $z$ is a directed edge, then the node $T(z)$
is called a {\bf child} of $S(z)$, while $S(z)$ is called a {\bf
parent} of $T(z)$.

A \textbf{subgraph} of a $\Gamma$-graph $G$ is a $\Gamma$-graph
$H$ each of whose nodes is a node in $G$ such that $H$ is closed
under children in $G$, and if $z$ is an edge in $G$ between nodes
in $H$, then $z$ is in $H$. A \textbf{family} is a subgraph
consisting of a node and all its children. A \textbf{path} through
the graph is a finite sequence of nodes $\langle n_0, \ldots , n_k
\rangle$, each connected to the next by a directed edge. If such a
path exists then $n_k$ is a \textbf{descendant} of $n_0$.  A {\bf
pointed $\Gamma$-graph} is a pair $(G,n^*)$ where $G$ is a
$\Gamma$-graph and $n^*$ is a distinguished node in $G$.

The nodes in our graphs will be  {\bf propositional letters} and
{\bf labelled $k$-ary boolean operators:} pairs $(p,\gamma)$ where
$p$ is a label (in some alphabet set) and $\gamma:\{ t,f\}^k\rightarrow \{t,f\}$
for some $k\geq 0$. (We can think of
propositional letters as labelled $-\infty$-ary boolean operators
so that all nodes are labelled operators.)

\begin{dfs}\label{formulaDefinition}
Let $\Gamma$ be the collection of all labelled boolean operators
and propositional letters. A {\bf formula} is a triple $\Phi =
(G,n^*, F)$ where:
\begin{enumerate}
\item $(G,n^*)$ is a pointed $\Gamma$-graph such that nodes with
$k$ children are $k$-ary boolean operators for $k>0$ and nodes
with no children are either constant boolean operators or
propositional letters.

\item $F$ is a set of nodes of $G$ which are said to be
\textbf{free}.
\end{enumerate}
In order that the interpretation of a graph not be ambiguous, we
assume that the set of children of each node is ordered. A pointed
subgraph of a formula $\Phi$ is called a {\bf subformula} if its free nodes are free in $\Phi$.
Its distinguished node need not coincide with that of $\Phi$.
\end{dfs}

\begin{remark}
In formulas of classical sentential logic, the propositional
letters can be thought of as free variables. We have generalized
this idea by designating arbitrary nodes in a formula as free.
\end{remark}

Note that each propositional letter can occur only once in a
formula, whereas the same boolean operator can occur
multiple times with different labels. Following are
two simple examples of formulas. In displaying a formula, we omit
the labels on the boolean operators, we arrange the arrows originating at each boolean operator from left to right to reflect the order of their children (when order is important, as in the case of non-symmetric
operators), we place a ``$*$" next to the distinguished node, and
we show free nodes in double circles:

\vskip .25 cm
\hskip 5.2cm $a$ \hskip 2.2cm \ $b$\\
\centerline { \xymatrix{& *+[o][F]{\rightarrow}\ar[dr]^(.1){\ *}\ar[dl]\\
*+[o][F=]{p}  && *+[o][F=]{q}}\qquad
 \xymatrix{
  *+[o][F]{\vee} \ar@(ur,dr)[]\ar[d]_(0){*\ }\\
*+[o][F=]{p} }
}

\begin{remarks}\hfill
\begin{enumerate}
\item Formula $a$ is the classical formula $p \rightarrow q$.

\item Cycles in the underlying graph permit us to encode
self-referential formulas. For instance, $b$ above can be read as
``$p \vee b$''.
 \item We can encode
formulas in classical sentential logic (such as $a$ above) using
trees as the underlying graphs (see below).
 \item Although we disallow duplicate copies of
 propositional letters in our graphs, classical formulas mentioning a
 propositional letter more than
  once can obviously still be realized by identifying the
   corresponding nodes.
\end{enumerate}
\end{remarks}

Barwise and Etchemendy construct a syntax for self-referential
statements in \cite{Barwise}. Following is a list of some of their
examples showing their notation and our equivalent representation
of these statements as formulas.

\begin{exs}\label{barwiseexamples}\hfill
\begin{enumerate}
\item The Liar (``This proposition is false.''):\\
$\phi =[ \mathit{Fa} \ \phi ]$

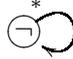
\centerline{\xymatrix{
 *+[o][F]{\neg} \ar@(ur,dr)[]^(.1){\ *}}}

 \vskip .25cm

\item The Strengthened Liar (``The Liar is false.''):

$\phi =[ \mathit{Fa} \ \phi ]$

$\psi =[ \mathit{Fa} \ \phi ]$

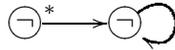
\centerline{\xymatrix{
 *+[o][F]{\neg}\ar[r]^(.2){\ *} & *+[o][F]{\neg} \ar@(ur,dr)[]}}
 \vskip .25cm

 \item Liar Cycle of Length Three (``The next proposition is true.
 The next proposition is true. The first proposition is
 false.''):\\
$\phi_1 =[ \mathit{Tr} \ \phi_2 ]$\\
$\phi_2 =[ \mathit{Tr} \ \psi ]$\\
$\psi =[ \mathit{Fa} \ \phi_1 ]$

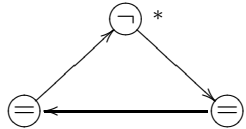
\centerline { \xymatrix{
& *+[o][F]{\neg}\ar[dr]^(.1){\ *}\\
*+[o][F]{=}\ar[ur]  && *+[o][F]{=}\ar[ll]}\qquad }

\vskip .25cm

\item The Contingent Liar (``Max has the three of clubs and this
proposition is false.''):\\
$\phi =[ \mathit{Max} \ H \ 3C] \wedge [ \mathit{Fa} \ \phi ]$

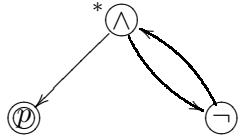
 \centerline {\xymatrix{
& *+[o][F]{\wedge}\ar[dl]_(0){*\ }\ar @/_/ [dr]\\
*+[o][F=]{p} && *+[o][F]{\neg}\ar@/_/ [ul]}
 }

\vskip .25 cm

\noindent Note that the propositional letter $p$ corresponds to
the proposition $[ \mathit{Max} \ H \ 3C]$.

\item Contingent Liar Cycle (``Max has the three of clubs. The
next proposition is true. At least one of the first two
propositions is false.''):\\
$\phi_1 =[ \mathit{Max} \ H \ 3C]$\\
$\phi_2 = [ \mathit{Tr} \ \psi ]$\\
$\psi = [ \mathit{Fa} \ \phi_1 ] \vee [ \mathit{Fa} \ \phi_2 ]$

 
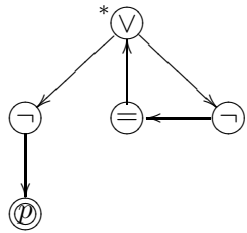
\centerline {\xymatrix{
& *+[o][F]{\vee}\ar[dl]_(0){*\ }\ar[dr]\\
*+[o][F]{\neg}\ar[d] & *+[o][F]{=}\ar[u] & *+[o][F]{\neg}\ar[l]\\
 *+[o][F=]{p} }
 }

\vskip .25 cm

\noindent Again, $p$ corresponds to $[ Max \ H \ 3C]$.

\vskip .25 cm

\item\label{lobsparadox} L$\ddot{\rm o}$b's Paradox (``If this
proposition is true,
then Max has the three of clubs.''):\\
$\phi =[ \mathit{Fa} \ \phi ] \vee [\mathit{Max} \ H \ 3C]$

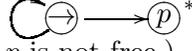
\centerline{\xymatrix{
 *+[o][F]{\rightarrow}\ar@(dl,ul)[]\ar[r]^(1.2){\ *} & *+[o][F]{p} }}

 \noindent (Note that the node $p$ is not free.)

 \vskip .25cm

\item Gupta's Puzzle (``Max has the three of clubs. The last two
propositions are true. At least one of the last two propositions
is false. Claire has the three of clubs. At most one of the first
three propositions is true.''):\\
$\phi_1 =[ \mathit{Max} \ H \ 3C]$\\
$\phi_2 = [ \mathit{Tr} \ \psi_1 ] \wedge [ \mathit{Tr} \ \psi_2 ]$\\
$\phi_3 = \overline{\phi_2}$\\
$\psi_1 =[ \mathit{Claire} \ H \ 3C]$\\
$\psi_2 = ([ \mathit{Fa} \ \phi_1 ] \wedge [ \mathit{Fa} \ \phi_2
]) \vee ([ \mathit{Fa} \ \phi_1 ] \wedge [ \mathit{Fa} \ \phi_3 ])
\vee ([ \mathit{Fa} \ \phi_2 ] \wedge [ \mathit{Fa} \ \phi_3 ])$

 
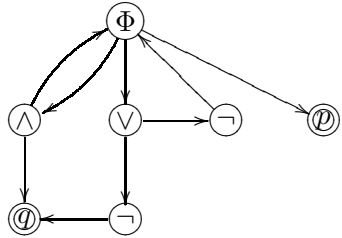
\centerline {\xymatrix{
& *+[o][F]{\Phi}\ar @/^/ [dl]\ar[drr]\ar[d]\\
*+[o][F]{\wedge}\ar @/ ^/ [ur]\ar[d] & *+[o][F]{\vee}\ar[r]\ar[d] & *+[o][F]{\neg}\ar[ul] & *+[o][F=]{p}\\
 *+[o][F=]{q} &  *+[o][F]{\neg}\ar [l]}
 }

\vskip .25cm

\noindent $p$ corresponds to the proposition $[ \mathit{Max} \ H \
3C]$ and $q$ to the proposition $[ \mathit{Claire} \ H \ 3C]$.
$\Phi$ is the boolean operator that says ``at most one of these
three is true." The location of the star will change depending
upon which truth value we wish to compute (see below).

\end{enumerate}
\end{exs}

\begin{dfs}\label{freenode} An \textbf{evaluation} of a formula $\Phi$ is a function $e$
from the set of free nodes to the set $\{T,F,V,L\}$. A node that
is mapped to $T$ or $F$ is said to be {\bf bound as true} or {\bf
false} respectively. A \textbf{proposition} is a pair $(\Phi , e)$
where $\Phi$ is a formula and $e$ is an evaluation.

If $(\Phi, e)$ is a proposition with underlying set of nodes $N$,
then a \textbf{hypothesis} on $(\Phi ,e)$ is a function
$H:N\rightarrow \{ t,f \}$. (Note that $H$ is not required to take
nodes bound as true (respectively false) to $t$ (respectively
$f$). However, we shall see shortly that it will rapidly evolve to
one that does.)

If $H$ is a hypothesis on $(\Phi , e)$ such that $H(n^*)=t$
(respectively $f$), then we say that the formula or proposition is
\textbf{assumed true} (respectively \textbf{assumed false}).
\end{dfs}

To describe rules for changing the evaluations on unbound nodes,
we need some notation. Let $c$ be a node with children
$\{\,c_1,\dots,c_n\,\}$ (so that $c$ is an $n$-ary boolean
operator
--- note that one of the $c_i$ may be $c$ itself) and denote by $\rho =
(\epsilon_{c_1},\dots ,\epsilon_{c_ n}; \epsilon_c)$ a row in the
truth table of $c$. Thus, each $\epsilon$ is either $t$ or $f$,
and $\epsilon_c$ is the {\bf output} value of $c$ determined by
the {\bf inputs} $\epsilon_{c_i}$. Note that permitting loops in
the graphs has the effect that rows of the truth table may include
truth values for the same node in two slots: as input and output.
If $C$ is a nonempty collection of children of $c$, denote by
$\rho|C$ the sub-tuple $(\epsilon_{c_i}; \epsilon_c)$ of $\rho$
with inputs indexed on the children in $C$. If $H$ is a hypothesis
on a proposition whose underlying graph includes $c$ and its
family, denote by $H|C$ the corresponding tuple $(\mu_{c_i};
\mu_c)$ of $t$s and $f$s determined by $H$.

\begin{df}\label{elemcon} Let $H$, $K$ be hypotheses on a proposition
$(\Phi , e)$. We say that $K$ is an \textbf{elementary
consequence} of $H$, and write $H\rhd K$, if $K$ is obtained from
$H$ by changing its value on a single node $c$ in one of the
following ways:
\begin{enumerate}

\item If $c$ is bound as $T$ and $H(c) = f$, then $K(c) = \neg
H(c)$.

\item If $c$ is bound as $F$ and $H(c) = t$, then $K(c) = \neg
H(c)$.

\item If $c$ is bound as $L$, then $K(c) = \neg H(c)$.

\item If $c$ is not bound and is a parent (possibly of itself),
let $C$ be a nonempty subset of its children, and let $\rho$ be a
row of the truth table of $c$ such that $H|C \not= \rho|C$, but
does agree if we change only the \emph{output} coordinate $\mu_c$
of $H|C$ to its negation. Then $K(c) = \neg H(c)$.

\item If $c$ is not bound and is a child (possibly of itself), let
$d$ be a parent, let $C$ be a subset of the children of $d$
containing $c$, and let $\rho$ be a row of the truth table of $c$
such that $H|C \not= \rho|C$, but does agree if we change only the
\emph{input} coordinate $\mu_c$ of $H|C$ to its negation. Then
$K(c) = \neg H(c)$.
\end{enumerate}
Note that no elementary consequence can change the value of any
node bound as $V$. We say that $K$ is a \textbf{consequence} of
$H$, or that \textbf{the assumption} $H$ \textbf{leads to the
conclusion} $K$, and write $H\rightarrowtail K$, iff $K$ is
obtained from $H$ by following a finite sequence of elementary
consequences.
\end{df}

\begin{exs}\hskip 2.7cm  \hskip 4.2cm \ \\
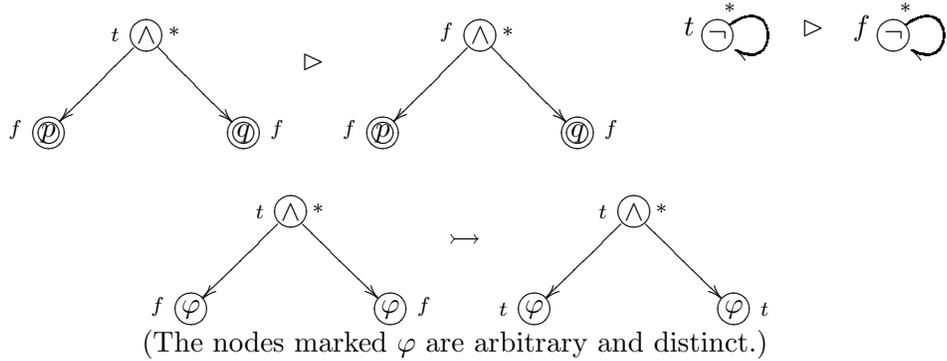
\centerline{
 \xymatrix{
&  *+[o][F]{\wedge}\ar[dr]^(.1){\,*}^(1.15){\,f}\ar[dl]_(.15){t\,}_(1.15){f\,}\\
*+[o][F=]{p}  && *+[o][F=]{q}} $\begin{matrix} \ \\ \ \\ \ \rhd
 \end{matrix}$  \xymatrix{
&  *+[o][F]{\wedge}\ar[dr]^(.1){\,*}^(1.15){\,f}\ar[dl]_(.15){f\,}_(1.15){f\,}\\
*+[o][F=]{p}  && *+[o][F=]{q}} \qquad
 $t\,$\xymatrix{
 *+[o][F]{\neg} \ar@(ur,dr)[]^(.1){\ *}\\}
 \qquad $\rhd$\ \ \
 $f\,$\xymatrix{*+[o][F]{\neg} \ar@(ur,dr)[]^(.1){\ *}}}

\centerline{ \xymatrix{
&  *+[o][F]{\wedge}\ar[dr]^(.1){\,*}^(1.15){\,f}\ar[dl]_(.15){t\,}_(1.15){f\,}\\
*+[o][F]{\varphi}  && *+[o][F]{\varphi}} $\begin{matrix} \ \\ \ \\
\ \rightarrowtail
 \end{matrix}$  \xymatrix{
&  *+[o][F]{\wedge}\ar[dr]^(.1){\,*}^(1.15){\,t}\ar[dl]_(.15){t\,}_(1.15){t\,}\\
*+[o][F]{\varphi}  && *+[o][F]{\varphi}}}
  \centerline{(The nodes marked $\varphi$ are arbitrary and distinct.)}
\end{exs}
\vskip .5 cm

Note that one of the following must hold for any proposition
$(\Phi , e)$:
\begin{enumerate}[ ]
\item \textbf{(T)} Every assumption that $(\Phi , e)$ is false
leads to the conclusion that it is true, and \textbf{not} every
assumption that $(\Phi , e)$ is true leads to the conclusion that
it is false.

\item \textbf{(F)} Every assumption that $(\Phi , e)$ is true
leads to the conclusion that it is false, and \textbf{not} every
assumption that $(\Phi , e)$ is false leads to the conclusion that
it is true.

\item \textbf{(V)} Not every assumption that $(\Phi , e)$ is true
leads to the conclusion that it is false, and not every assumption
that $(\Phi , e)$ is false leads to the conclusion that it is
true.

\item \textbf{(L)} Every assumption that $(\Phi , e)$ is false
leads to the conclusion that it is true, and  every assumption
that $(\Phi , e)$ is true leads to the conclusion that it is
false.
\end{enumerate}

\begin{remarks}\hfill
\begin{enumerate}
\item In the cases that \textbf{not} every assumption that $(\Phi
,e)$ is true (respectively false) leads to the conclusion that it
is false (respectively true), we will sometimes say that $(\Phi
,e)$ can {\bf get stuck} in true (respectively false).

\item It follows from the definitions that every node bound as $T$
becomes $t$ under an elementary consequence, and its value is
subsequently fixed as $t$. Similarly, nodes bound as $F$ turn $f$
and remain so.
\end{enumerate}
\end{remarks}

The above observation leads us to a four-valued system of logic:

\begin{df}\label{truthvalue} A proposition $(\Phi , e)$ has a
\textbf{truth value} $P\in \{ T,F,V,L\}$ depending on which of the
above possibilities holds.

\end{df}

\begin{remarks} \hfill
\begin{enumerate}

\item $T$ stands for ``true", $F$ for ``false", $V$ for
``vacuous", and $L$ for ``lie."

\item The Liar (``This statement is false.'')

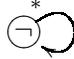
\centerline { \xymatrix{
 *+[o][F]{\neg} \ar@(ur,dr)[]^(.1){\ *}\\}\qquad
 }
\vspace{.5cm} \noindent has truth value $L$. On the other hand,
the \emph{Vacuous Affirmation} (``This statement is true."):

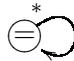
\centerline { \xymatrix{
 *+[o][F]{=} \ar@(ur,dr)[]^(.1){\ *}\\}\qquad
 }
\vspace{.5cm} \noindent has truth value $V$.

\item Let's compute the truth values of the other formulas in
Example \ref{barwiseexamples}. Note that they agree with the
heuristic reasoning found in Barwise and Etchemendy
\cite{Barwise}. The strengthened liar has truth value $L$, as does
any liar cycle. The contingent liar and contingent liar cycle have
truth value $F$ when the propositional letter is bound as $F$, and
truth value $L$ when the propositional letter is bound as $T$.
L$\ddot{\rm o}$b's Paradox has truth value $T$ (Note that in this
example no nodes are free and we are computing the truth value of
the propositional letter $p$). In Gupta's puzzle, under the
assumption that Claire has the three of clubs and Max does not we
bind $p$ as $F$ and $q$ as $T$. By starring each node in
succession we see that the nodes corresponding to $\wedge$ and
$\Phi$ have truth value $T$, and the node corresponding to $\vee$
has truth value $F$.

\end{enumerate}

\end{remarks}

\begin{df} A formula is
\textbf{well-grounded} if no node is a descendant of itself and
every node is a descendant of the distinguished node. It is
\textbf{strongly well-grounded} if, in addition, the free nodes
are precisely the propositional letters.
\end{df}

Note that there is a natural bijection $\Psi :W\rightarrow P$ from
the collection $W$ of strongly well-grounded formulas to the
collection $P$ of formulas in sentential calculus. If $e$ is an
evaluation of $\Phi \in W$ with the propositional letters bound as
$T$ or $F$, then $\Psi (\Phi )$ has a truth value coinciding with
the truth value of $(\Phi , e)$ under the assignment of values to
the propositional letters under $\Psi$ (note that the truth value
under these circumstances is always $T$ or $F$). We can therefore
identify the strongly well-grounded formulas with the
corresponding formulas of sentential calculus.

\begin{ex} $p\leftrightarrow \neg p$ can be identified with the
following well-grounded formula:

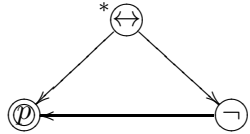
\centerline {\xymatrix{
& *+[o][F]{\leftrightarrow}\ar[dl]_(0){*\ }\ar[dr]\\
*+[o][F=]{p}  && *+[o][F]{\neg}\ar[ll]}\qquad
 }

\end{ex}

\section{Logical Equivalence and Classification of Propositions}\label{Logical Equivalence and Classification of Propositions
}

\begin{df} A formula $\Phi$
is \textbf{completely connected} if each boolean operator in
$\Phi$ has every node (including itself) as a child.
\end{df}

Given a formula $\Phi$ with $k$ nodes, we can define a completely
connected formula $\Phi^{\prime}$ by replacing each boolean
operator $\gamma: \{t, f\}^q \to \{t, f\}$ with a $k$-ary boolean
operator $\gamma^{\prime} :\{t, f\}^k \to \{t, f\}$ obtained by
composing $\gamma$ with the evident projection $\pi :\{t, f\}^k
\to \{t, f\}^q$. Note that $\Phi$ and $\Phi^{\prime}$ are
indistinguishable in the sense that corresponding hypotheses on
$\Phi$ and $\Phi^{\prime}$ have corresponding elementary
consequences

\begin{dfs} The \textbf{truth table} of a formula with $k$ free
nodes is the function $h:\{ T,F,L,V \}^{k}\rightarrow \{ T,F,L,V
\}$ obtained by computing the truth value of the starred node
resulting from each evaluation of the free nodes. (Implicit here
is an ordering of the free nodes, since we are thinking of $\{T,F,L,V \}^{k}$
as the collection of all evaluations of the free
nodes.) If a formula contains no free nodes then its truth table
consists of its (single) truth value. The \textbf{restricted truth
table} of a formula  is the restriction of its truth table to
$\{T, F\}^{k}$.

\end{dfs}

\begin{exs}\hfill
\begin{enumerate}
\item The truth table of the Liar is $L$.

\item If $\Phi$ is any formula with $k$ free nodes in which the
distinguished node $n^*$ is free, then its truth table is the
projection $\pi :\{ T,F,L,V \}^{k}\rightarrow \{ T,F,L,V \}$ onto
the coordinate associated with $n^*$.

\item Following is a representation of the truth table $\{ T,F,L,V
\}^2\rightarrow \{ T,F,L,V \}$ of the formula $p\wedge q$:

\vskip .25cm  \centerline{
\begin{tabular}{|c|cccc|}
  \hline
  \  & T & F & L & V \\
  \hline
  T & T & F & L & V \\
  F & F & F & F & F \\
  L & L & F & L & F \\
  V & V & F & F & V \\
  \hline
\end{tabular}}
 \vskip .25 cm

\end{enumerate}
\end{exs}

\begin{dfs}\label{logicallyequivalent} Two formulas with the same free nodes are \textbf{logically
equivalent} if they have the same truth table. They are
\textbf{weakly logically equivalent} if they have the same
restricted truth table. A {\bf (strong) tautology} is a formula
whose truth table has constant value $T$ while a {\bf weak
tautology} is a formula whose restricted truth table has constant
value $T$
\end{dfs}

The above definitions beg the following questions:
\begin{enumerate}
\item Are weakly logically equivalent formulas logically
equivalent?

\item Are all truth tables realizable? That is, is every function
$$h:\{ T,F,L,V \}^{k}\rightarrow \{ T,F,L,V \}$$ the truth table
of some formula with $k$ free nodes?

\item Are all restricted truth tables realizable?

\item If $\Psi$ is a subformula of $\Phi$, and we replace $\Psi$
by any equivalent subformula $\Psi'$, is the resulting formula
logically equivalent to $\Phi$?
\end{enumerate}

The easiest question to answer is the first. The formulas
\begin{equation*}
p \vee \neg p \qquad \text{ and } \qquad  p \leftrightarrow p
\end{equation*}
are weakly logically equivalent but not logically equivalent; the
first returns $L$ when $p$ is bound as $L$, while the second
returns $T$.

We shall answer the second question negatively below, and in the
process describe which truth tables are realizable. To do so
involves first answering the third question affirmatively.

\begin{prop}\label{realization1} All restricted truth tables are realizable. That is,
every function $h:\{ T,F\}^{k}\rightarrow \{ T,F,L,V\}$ is the
restricted truth table of some formula $\Phi$.
\end{prop}

\begin{proof}
The desired formula $\Phi$ will be completely connected with $k$
free nodes $n_1,\dots,n_{k}$ that are propositional letters and a
single node $n^*$ which is not free. The distinguished node $n^*$
is a (labelled) $(k+1)$-ary boolean operator $\varphi$, specified
by
\begin{equation*}
\varphi(\epsilon_1,\dots,\epsilon_{k},\epsilon) =
\begin{cases} t & \text{ if } h(\epsilon_1,\dots,\epsilon_{k}) = T\\
f & \text{ if }h(\epsilon_1,\dots,\epsilon_{k}) = F\\
\epsilon & \text{ if }h(\epsilon_1,\dots,\epsilon_{k}) = V\\
\neg\epsilon & \text{ if }h(\epsilon_1,\dots,\epsilon_{k}) = L
\end{cases}.
\end{equation*}
Here the last coordinate corresponds to the starred node. It is
then easy to check that $\Phi$ behaves as desired.
\end{proof}

Next, we turn to the question of which truth tables can be
realized. First, partially order the four truth values as follows:

\begin{df} $\prec$ is the partial ordering defined
on $\{T,F,L,V \}$ by the following diagram:
$$V\prec F, T\prec L$$
\end{df}
Thus, formulas with lesser truth values are more prone to getting
stuck.

 If $(\epsilon_1,\dots, \epsilon_{k}) \in
\{T,F,L,V\}^{k}$, then let $\mathcal L$ (resp. $\mathcal V$) be
the set of indices $i$ for which $\epsilon_i = L$ (resp. $V$).
First, the definition of truth leads to the following observation:

\begin{prop}\label{inequality}
If $\Phi$ is a formula with $k$ free nodes and restricted truth
table $h\colon \{T,F\}^{k} \to \{T,F,L,V\}$ and if $(\epsilon_i)
\in \{T,F,L,V\}^{k}$, then
\begin{equation*}
\begin{split}
&h((\epsilon_i)) \geq \\ &\inf_{\mu_i \in \{T ,F\}, i \in \mathcal
V}\ \sup_{\mu_i \in \{T ,F\}, i \in \mathcal L}   \{h(\mu_1,\dots,
\mu_{k})\mid \mu_i = \epsilon_i \text{ if } \epsilon_i \in \{T
,F\}\}
\end{split}
\end{equation*}
\end{prop}
\begin{proof}
From the definition of truth, replacing one or more bound $T$ or
$F$ nodes by bound $L$ nodes replaces the truth value of $\Phi$ by
values at least as large as the supremum of the original values.
Subsequently replacing one or more of the remaining nodes by $V$
nodes results in exactly the infimum of the current values.
\end{proof}

\begin{remarks}\hfill
\begin{enumerate}
\item Note that reversing $\sup$ and $\inf$ in the above
proposition gives a lower bound, and hence a weaker result.

\item The construction in Proposition \ref{realization1} produces
formulas whose values on elements of $\{T,F,L,V\}^{k}$ can be
computed by replacing the inequality above by an equality. Thus,
it must be shown that we can introduce additional appropriate
elementary consequences in the construction there without
affecting the values on $\{T,F\}^{k}$.
\end{enumerate}
\end{remarks}

\begin{thm}\label{realization}
Every instance of the inequality in Proposition \ref{inequality}
can be realized.
\end{thm}
\begin{proof}
Suppose $h$ satisfies the inequality in Proposition
\ref{inequality}, let $g$ be the corresponding restricted truth
table, and let $\Phi$ be the realization of $g$ constructed in
Proposition \ref{realization1}. It suffices to show that we can
adjust $\Phi$, without altering its restricted truth table, in
such a way that the starred node can change in a prescribed way
(from $t$ to $f$, from $f$ to $t$ or both) when any specified
collection $S$ of free nodes is bound as $L$. Denote this desired
elementary consequence (on the starred node) by $c_S$.

We accomplish this by expanding $\Phi$ as follows. First, add a
new free node $n_S$ for every nonempty subset $S$ of nodes in the
original collection of free nodes. Next, construct a new formula
$\Psi$ whose restricted truth table agrees with that of $\Phi$ on
all evaluations in which all the new nodes are set to $f$, but
causes the desired change $c_S$ in the starred node precisely when
$n_S$ is $t$ and $n_{S'} = f$ for all $S' \neq S$.

Next, we turn each new free node in $\Psi$ into a ``lie detector''
as follows: for each $S$, write $S = \{n_1, \dots, n_r\}$ and
replace $n_S$ by
\begin{equation}\label{new}
[n_1\wedge \dots \wedge n_r] \wedge \neg[n_1\vee \dots \vee n_r]
\end{equation}
Thus, the new nodes are no longer free. Call the resulting formula
$\Psi'$. If the original free nodes are bound as $T$ or $F$, and
if $H$ is any hypothesis on $\Phi$ that causes it to become stuck,
then that same hypothesis on $\Psi'$ with each $n_S$ set to false,
and the conjunctions and disjunctions in (\ref{new}) set to the
appropriate values, is still stuck. Therefore, $\Phi$ and $\Psi'$
have the same restricted truth table. On the other hand, if any
subset $S$ of free nodes is bound as $L$, then $n_S$---and only
$n_S$--- can turn to $t$ in any configuration in which the old
circuit was stuck. But since turning $n_S$ on has the desired
effect, we are done.
\end{proof}

\begin{ex}It is instructive to consider the special case $k = 1$.
Call a formula with one free node that is also a propositional
letter a {\bf gate.} More specifically, if $P_1, P_2, P_3$ and
$P_4$ are truth values, then a $P_1P_2P_3P_4${\bf -gate} is a gate
with truth table $T\mapsto P_1, F\mapsto P_2, L\mapsto P_3,
V\mapsto P_4$. For instance, negation can be viewed as an
$FTLV$-gate. \cite{Scott} has shown that there are exactly 25
gates out of a possible 256, and has constructed simple formulas
of each type. One non-existent gate is an $FTVL$-gate.

\end{ex}

 To end this section, we turn to the question of replacing
subformulas by equivalent ones. Since it is possible for a formula
to have edges pointing to different nodes in a subformula, we
shall restrict the kinds of subformulas we can swap.

\begin{df} If $\Phi$ and $\Psi$ are any two formulas with no nodes in common,
and $p$ is a propositional letter in $\Phi$, then the formula
obtained by replacing the node $p$ by the formula $\Psi$ is called
the \textbf{substitution} of $\Psi$ into $\Phi$ at node $p$, and
is written as $s(\Psi , \Phi, p)$
\end{df}

\begin{prop} Let $\Phi$, $\Psi_1$, $\Psi_2$ be
formulas such that $\Phi$ has no nodes in common with $\Psi_1$ or
$\Psi_2$, and $\Psi_1$ and $\Psi_2$ are logically equivalent. Let
$p$ be a propositional letter in $\Phi$. Then $\Phi_1=s(\Psi_1 ,
\Phi , p)$ and $\Phi_2=s(\Psi_2 , \Phi , p)$ are logically
equivalent.
\end{prop}

\begin{proof}Let $e$ be an evaluation of $\Phi_1$ and
suppose that the proposition $(\Phi_1,e)$ has truth value $P$.
That $(\Phi_2,e)$ also has truth value $P$ follows easily from
Definition ~\ref{elemcon}. We show one case as an illustration:
Suppose that $\Psi_1$ has truth value $L$ under the current
evaluation (thus $\Psi_2$ also has truth value $L$), and suppose
that $(\Phi_1,e)$ can get stuck in $T$. Let $H_1$ be a hypothesis
on $(\Phi_1,e)$ such that $H_1(n^*)=t$ and there is no consequence
$H_1^{\prime}\leftarrowtail H_1$ with $H_1^{\prime}(n^*)=f$.
Define the hypothesis $H_2$ on $(\Phi_2,e)$ by $H_2(x)=H_1(x)$ for
all nodes $x\ne p$ in $\Phi$  and let $H_2(x)$ be arbitrary for
$x$ in $\Psi_2$. It is then easy to see that there is no
$H_2^{\prime}\leftarrowtail H_2$ with $H_2^{\prime}(n^*)=f$. The
other cases are just as easy.
\end{proof}

\section{Models of Sentential Logic}\label{Models of Sentential Logic}

It may seem most natural to generalize the notion of a model (see,
for example, \cite{sententialRef}) by defining a model to be an
assignment  of truth values to the propositional letters. However,
we wish to create models which admit inconsistencies of an
arbitrary nature, for instance models in which a propositional
letter and its negation are both true. We would also like to
create models in which natural inconsistencies such as the Liar
can become true or false. The following definition serves our
needs in this paper, and will also serve as the basis for our
mathematical model theory in the sequel \cite{Part2}.

\begin{df} A \textbf{model} $\mathcal M$ of sentential logic
is a tuple $(M,\mathcal A, \mathcal B,\mathcal C, \mathcal D)$
where $M$ is a nonempty set of propositional letters, and
$\mathcal A, \mathcal B$, $\mathcal C$ and $\mathcal D$ are
disjoint sets of formulas without free nodes.
\end{df}

\begin{remarks}\hfill
\begin{enumerate}

\item The formulas in $\mathcal A$ are called \textbf{true
axioms}. Similarly, the formulas in $\mathcal B, \mathcal C$ and
$\mathcal D$ are called \textbf{false axioms, contradictory axioms
and independent axioms}, respectively. When determining truth,
each axiom will be bound to have the truth value determined by
which of these four sets it is a member (see below).

\item Since the nodes in our formulas are \emph{labelled} boolean
operators, it follows that models are possible in which different
instances (distinguished by their labels) of the same formula are
assigned different truth values as axioms. For instance, some
Liars could be true and others false.

\end{enumerate}
\end{remarks}

\begin{df}
Let $\mathcal M$ be a model. A {\bf formula} $\Phi$ {\bf in}
$\mathcal M$ is a formula with no free nodes whose propositional
letters are all in $M$.
\end{df}

We now show how to compute the truth value of a formula $\Phi$ in
$\mathcal M$.

\begin{df}
A {\bf subaxiom} of an axiom $A$ in $\mathcal M$ is a proper
subformula $B$ of $A$.
\end{df}

\begin{df}\label{truthdef}The truth values of formulas in $\mathcal M$
are determined as follows:
\begin{enumerate}

\item If $\Phi$ is an axiom, then its truth value is determined by
its membership in $\mathcal A, \mathcal B$, $\mathcal C$ or
$\mathcal D$ as above.

\item If $\Phi$ is a subaxiom but not an axiom, its truth value
$R(\Phi )$ will be defined as a limit of a nondecreasing sequence
of truth values $R_n(\Phi )$, as follows: For any axiom $A$ of
which $\Phi$ is a subaxiom let $A^*$ be $A$ with the star moved to
the node for $\Phi$, bind all axioms by their truth values under
(1) that occur as subformulas of $A$, and compute the truth value
of $A^*$. Denote this truth value by $R_A$. Then take
\begin{equation*}
R_1(\Phi ) = \sup \{ R_A \mid \Phi\text{  a subformula of  }A\}.
\end{equation*}
To obtain $R_n(\Phi )$ from $R_{n-1}(\Phi )$, proceed as in the
definition of $R_1$, but also bind all subaxioms other than $\Phi$
and occurring as subformulas of $A^*$ by their truth value under
$R_{n-1}$. Let $R(\Phi )=\sup_{n<\omega}R_n(\Phi )$.

\item If $\Phi$ is not an axiom or subaxiom, bind all its
subformulas that are axioms or subaxioms by their truth values
under (1) and (2), and then compute the value of the starred node
as usual.
\end{enumerate}

\begin{remarks}
\begin{enumerate}

\item Technically, in computing truth we are using formulas
identical to the ones we are interested in, but with certain nodes
free so that we can bind them.

\item Note that we do not permit formulas that are not subaxioms
to inherit truth values from formulas containing them as we do in
Definition \ref{truthdef}(2). If we did then the presence of a
single $L$ formula would have the consequence that every formula
that is not a subaxiom would have truth value $L$. Indeed, if
$\Phi$ is $L$ and not a subaxiom, and $\Psi$ is any other formula
that cannot get stuck in $F$, then $\Phi \wedge \Psi$ is seen to
be $L$. But binding this formula as $L$ leads to $\Psi$ having
truth value $L$. If $\Psi$ can get stuck in $F$ then a similar
argument works using a disjunction instead of a conjunction. By
contrast Definition \ref{truthdef} will allow us to admit
contradictions and contain them.

\end{enumerate}
\end{remarks}

 We write
$\mathcal M \models \Phi$, $\mathcal M \models^F \Phi$, $\mathcal
M \models^L \Phi$ or $\mathcal M \models^V \Phi$ if $\Phi$ takes
on truth value $T, F, L$ or $V$ respectively. $\mathcal M$ is
\textbf{complete} if there are no well-grounded formulas $\Phi$ in
$\mathcal M$ such that $\mathcal M\models^V \Phi$ (in particular,
$\mathcal D$ contains no well-grounded formulas).
\end{df}
\begin{ex}

Let $\mathcal M=(M,\mathcal A, \mathcal B,\mathcal C, \mathcal D)$
where $M=\{ p\}$, $\mathcal A=\{ p,\neg[=p]\}$ and $\mathcal
B=\mathcal C=\mathcal D=\emptyset$. Then the subaxiom $[=p]$ has
truth value $L$ in $\mathcal M$ as the following diagram shows.

\centerline {\xymatrix{  *+[o][F]{\neg}\ar[r]_(.3){T\ } &
*+[o][F]{=}\ar[r]_(.3){*\ }_(.8){T\ } & *+[o][F]{p}}
 }

\end{ex}

We now relate our models to classical models of sentential
calculus.

\begin{dfs}Let $\Phi$ be any formula in the model $\mathcal M$. Say
that $\Phi$ is \textbf{generically inconsistent} if its truth
value in $\mathcal M$ differs from its truth value obtained by
binding only the propositional letters that occur as axioms by
their truth values in $\mathcal M$.

Also, call a well-grounded formula $\Phi$ in the model $\mathcal
M$ {\bf special in} $\mathcal M$ if each node in $\Phi$ is either
a conjunction, disjunction, negation, equals, or a propositional
letter, and such that no subformulas of $\Phi$ of the form
``$=\Psi$'' are axioms. Further, we require that, if $\Psi$ is a
subformula of $\Phi$ that is a subaxiom of some axiom $A$, then
the formula ``$=\Psi$'' is also a subaxiom of $A$ (so that the
``$=$'' node is bound to at least the same value as its target
when we evaluate truth). In particular, the starred node of any
axiom in $\Phi$ must point to only ``$=$'' nodes.
\end{dfs}

\begin{lemma}\label{treelemma2}Let $\mathcal M$ be a model all of
whose axioms are either propositional letters assumed true or
false, or other well-grounded formulas assumed true. Furthermore,
assume that every propositional letter in $M$ is an axiom. If
there is a generically inconsistent special formula $\Phi$ in
$\mathcal M$, then $\Phi$ has an $L$ subformula.
\end{lemma}

\begin{proof}
Note that, since each propositional letter in $M$ is assigned a
unique truth value, $\mathcal M$ determines an evaluation $e$ of
every formula whose propositional letters are free and in
$\mathcal M$. Let $\Phi$ be as in the hypothesis and assume that
$\Phi$ has no $L$ subformulas. It follows that the truth value of
any subaxiom which is a subformula of $\Phi$ can be computed by
choosing any axiom of which it is a subformula, binding all
axioms, and then computing the truth value of the subaxiom's
distinguished node as usual.

Let $\Psi$ be a minimal subformula of $\Phi$ which is generically
inconsistent under the evaluation $e$ determined by $\mathcal M$.
Then $\Psi$ must be an axiom or subaxiom by minimality. If the
distinguished node $n^*$ of $\Psi$ is an ``$=$'' node, then the
definition of a special formula implies that it either has the
same truth value as its target, contradicting minimality, or is
$L$, contrary to assumption. Thus $n^*$ is not an ``$=$'' node.

Since $\Psi$ must be an axiom or a subaxiom, the special property
of $\Phi$ now guarantees that $n^*$ points to only ``$=$'' nodes,
each of which is the distinguished node of a subaxiom so that all
of these ``$=$'' nodes are bound (when we evaluate truth). Since
$\Psi$ is generically inconsistent, binding only the propositional
letters will lead any hypothesis to one in which $n^*$ has the
opposite truth value. Further the chain of elementary consequences
can be arranged to affect $n^*$ only in the last step (because
$\Psi$ is well-grounded).

If $n^*$ and the propositional letters are bound, then the same
chain of elementary consequences except for the last is still
possible. Since proper subformulas of $\Psi$ are generically
consistent, all nodes except $n^*$ end up in their values as
subaxioms. Therefore, binding $n^*$ to its given truth value will
result in an elementary consequence in which the truth value of
one of its children $c$ changes to the opposite truth value of its
associated subaxiom $C$ (recall that $n^*$ must be a conjunction,
disjunction or negation; see remark below).

One has $C < \Psi <A$ for some axiom $A$ (where $<$ indicates
subformula). A lower bound of the truth value of $C$ is obtained
by binding all axioms that occur in $A$ and then computing the
truth value of the node $c$. However, any hypothesis on $A$ will
lead to one in which $n^*$ is assigned its given truth value.
Subsequently, as we have seen, a further sequence of elementary
consequences will change $c$ to the opposite truth value. The
definition of truth value of subformulas now tells us that $C$
must be $L$ in the model, contrary to assumption.

\end{proof}

\begin{remark}

The definition of special need not be so restrictive as to include
only conjunctions, disjunctions, negations and equals. However, we
need to avoid operators which either are, or can behave like,
constant boolean operators. For example the operator
$\leftrightarrow$ can behave as a constant unary operator as in
$p\leftrightarrow p$. If this formula is hypothesized as false,
then there is no elementary consequence that changes the value of
its only child $p$.

\end{remark}

If $T$ is a set of well-grounded formulas, let $T'$ be obtained
from $T$ by writing each formula in disjunctive normal form, and
then replacing each non-starred node $n$ in each formula by
``$=n$''. Call this {\bf expanded disjunctive normal form.}
If $T'$ is now the set of axioms in a model
$\mathcal M$, then each axiom is special. It follows that any
classical theory $T$ is equivalent (in the sense that it has the
same models) to one with a special set of axioms.

\begin{thm}\label{classicalmodeltheorem} Let $T$ be a theory and
$M$ a model in the classical sense. Form the model $\mathcal M(T)$
by letting $\mathcal A$ consist of the formulas of $T$ written in
expanded disjunctive normal form, together with the true
propositional letters in $M$, and $\mathcal B$ the false
propositional letters in $M$. Then $M$ is a model of $T$ in the
classical sense iff there are no well-grounded $L$ formulas in
$\mathcal M(T)$. When this is the case, the well-grounded formulas
in $\mathcal M(T)$ have the same truth values as the corresponding
formulas in $M$.
\end{thm}

\begin{proof}
First observe that if $T'$ is the theory $T$, but with all axioms
written in expanded disjunctive normal form, then $M$ is a model
of $T$ iff $M$ is a model of $T'$. Therefore we assume without
loss of generality that the formulas in $T$ are in expanded
disjunctive normal form.

If $M$ is a model of $T$, then it is clear that all well-grounded
formulas in $\mathcal M(T)$ get stuck in their appropriate truth
values ($T$ or $F$), and hence there are no well-grounded $L$
formulas. Conversely, if $M$ is not a model of $T$, then there
exists a well-grounded formula $\Phi$ such that $T\vdash \Phi$ but
$\Phi$ is false in the model $M$. This implies that at least one
axiom in $T$ is generically inconsistent as a formula in $\mathcal
M(T)$. Since each axiom is special, the lemma implies that this
axiom has an $L$ subformula.
\end{proof}

In \cite{Barwise} Barwise and Etchemendy describe a class of
models which include false liars. It is therefore natural to ask
if there is an interesting class of models in our sense in which
the Liar is false. Uninteresting models can be constructed by
simply binding every formula with an arbitrary truth value. In an
interesting model, truth values must, for the most part, be
computed according to Definition \ref{truthdef}.

\begin{df} If $\Phi_1,\dots, \Phi_n$ are (not necessarily well-grounded) formulas
and $\Psi$ is a well-grounded formula with $n$ propositional
letters, let $\Psi'$ be obtained from $\Psi$ by replacing its
propositional letters with the $\Phi_i$. Then we say that the pair
$(\Psi',(\Phi_i))$ is {\bf relatively well-grounded.}
\end{df}

\begin{prop} If $T$ is any consistent theory with a classical model $M$,
then it has a model $\mathcal N(T)$ in our sense in which all
well-grounded formulas have the same truth values as the
corresponding formulas in $M$. Further, there are no $L$ formulas
in $\mathcal N(T)$, the Liar is false in $\mathcal N(T)$, and no
relatively well-grounded formula is an axiom or subaxiom.
(Thus, for instance, the negation of the Liar is true).
\end{prop}

\begin{proof}
We construct the model $\mathcal N(T)$ as the $\omega$-limit of an
inductively defined sequence $\mathcal M_i$. Define $\mathcal M_0
= \mathcal M(T)$, as in Theorem \ref{classicalmodeltheorem}, and
assume that $\mathcal M_i = (\mathcal A_i, \mathcal B_i, \mathcal
C_i, \mathcal D_i),\ (i < n)$ have been defined, with $\mathcal
A_i = \mathcal A_{i+1}$, $\mathcal B_i \subseteq \mathcal
B_{i+1}$, and $\mathcal C_i = \mathcal D_i = \emptyset$, and such
that all well-grounded formulas have the same truth values as the
corresponding formulas in $M$. To define $\mathcal M_{i+1}$, take
$\mathcal A_{i+1} = \mathcal A_i,\ \mathcal C_{i+1} = \mathcal
D_{i+1} = \emptyset$, and
\begin{equation*}
\mathcal B_{i+1} = \mathcal B_i \cup \{\,\Phi \mid \Phi \text{  a
minimal  }L\text{  formula in  }\mathcal M_i\,\}
\end{equation*}
Since we are decreasing some truth values in passing from
$\mathcal M_i$ to $\mathcal M_{i+1}$, it follows that no new $L$
formulas are introduced. Further, since by assumption, no
well-grounded formula can be $L$, their truth values are not
affected. Also, if a relatively well-grounded formula is $L$ in
$\mathcal M_i$, it cannot be a minimal $L$ formula, and so it will
remain unbound when we evaluate truth. Finally, it is not hard to
see that there are no $L$ formulas in the limit.
\end{proof}

\begin{remarks}\hfill
\begin{enumerate}
\item We cannot eliminate all $L$ formulas in just one step. For
instance, the following $L$ formula remains $L$ in $\mathcal M_1$
because only the bottom left node is bound as $F$ in $\mathcal
M_1$.

 
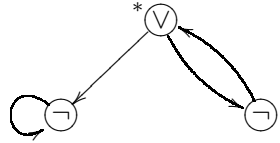
\centerline {\xymatrix{
& *+[o][F]{\vee}\ar[dl]_(0){*\ }\ar @/_/ [dr]\\
*+[o][F]{\neg}\ar@(ul,dl)[]  && *+[o][F]{\neg}\ar@/_/ [ul]}
 }
 \vskip .25 cm
\item We cannot avoid creating $V$ statements in general.  The
following formula has truth value $T$ until the Liar subformula
(rightmost node) is bound as $F$, whereupon the formula becomes
$V$. (Note also that we can replace the Liar on the right by any
$L$ formula.)

 
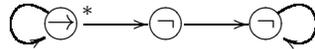
\centerline {\xymatrix{*+[o][F]{\rightarrow}\ar@(ul,dl)[]\ar[r]^(0.3){*\ }  &
 *+[o][F]{\neg}\ar[r] & *+[o][F]{\neg}\ar@(ur,dr)[]}}

 \vskip .25 cm
\item We cannot expect to eliminate $V$ statements by binding them
as $T$ or $F$ without (re-)introducing $L$ statements. In the
following pair of formulas, $\Phi$ stands for any $V$ subformula.
Binding $\Phi$ as false causes the first formula to change from
$T$ to $L$. On the other hand, binding $\Phi$ as true causes the
second formula to change from $T$ to $L$. Note that re-binding
these new $L$ statements will lead again to further $V$ statements
(parenthetical comment in (2)).

 
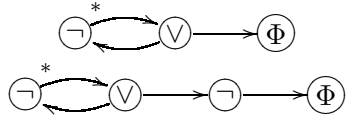
\centerline {\xymatrix{*+[o][F]{\neg}\ar @/^/ [r]^(0.3){*\ }  &
 *+[o][F]{\vee}\ar @/^/ [l]\ar[r] & *+[o][F]{\Phi}}}

 \centerline {\xymatrix{*+[o][F]{\neg}\ar @/^/ [r]^(0.3){*\ }  &
 *+[o][F]{\vee}\ar @/^/ [l]\ar[r] & *+[o][F]{\neg}\ar[r]&  *+[o][F]{\Phi}}}

 \vskip .25 cm

 \item Barwise  and Etchemendy \cite{Barwise} seem to avoid V
statements in their models of Austinian logic by limiting the
amount of information that is talked about in their
``situations.'' Our models, on the other hand, include all
formulas.

\end{enumerate}
\end{remarks}

We end this section with a brief discussion of rules of inference
in the hope of stimulating further research.

 Since
an instance of modus ponens may be bound as $F, V$, or $L$ in a
model, we cannot expect a general set of rules of inference to
hold in every model of sentential logic. However, models in which
all the axioms are propositional letters and bound as either $T$
or $F$ are interesting in this regard, since a large class of
rules of inference hold ``from without'' in such models.

\begin{df} Say a model is \textbf{simple} if each axiom is a
propositional letter assumed true or false, and every
propositional letter in $M$ is an axiom.
\end{df}

As usual, say that modus ponens \textbf{holds externally} if
$\mathcal M \models \Phi$ and $\mathcal M \models [\Phi
\rightarrow \Psi]$ implies $\mathcal M \models \Psi$, and
similarly for the other common rules of inference from classical
sentential logic. If $\Phi$ is a (weak or strong) tautology (see
Definition \ref{logicallyequivalent}), write $R(\Phi)$ for the
rule of inference that says $\mathcal M \models \Phi$.

\begin{prop} If $\mathcal M$ is a simple model, then the
following rules of inference hold: modus ponens, modus tollens,
contrapositive, chain rule, disjunctive inference, double
negation, De Morgan, simplification, conjunction, and disjunctive
syllogism. Further, the rule $R(\Phi)$ holds iff $\Phi$ is a
strong tautology.
\end{prop}
\begin{proof}The last statement is a consequence of our rules for computing
truth. In the case of modus ponens, suppose that $\mathcal M
\models \Phi$ and $\mathcal M \models [\Phi \rightarrow \Psi]$,
but that $\Psi$ has a truth value other than $T$. Then consider
the three possibilities: If $\mathcal M \models^F \Psi$, then
$\Phi \rightarrow \Psi$ has truth value $F$. If $\mathcal M
\models^L \Psi$, then $\Phi \rightarrow \Psi$ has truth value $L$.
Finally, if $\mathcal M \models^V \Psi$, then $\Phi \rightarrow
\Psi$ has truth value $V$. We therefore rule these possibilities
out by assumption. The remaining rules listed can be checked
one-by-one.
\end{proof}

Note that the modus ponens formula $[p \wedge (p\rightarrow
q)]\rightarrow q]$ is a (strictly) weak tautology, and so cannot
hold ``internally'' in any classical model. Further, not every
weak tautological implication leads to an external rule of
inference, as illustrated by the weak tautological implication
$(q\rightarrow(p\vee\neg p)]$. Indeed, if $\mathcal M \models q$,
then $\mathcal M \nvDash (p\vee\neg p)$ if $p$ is $L$, since then
$\mathcal M \models^L (p\vee\neg p)$.

In general, a good set of rules of inference and a good theory of
argument should include rules that predict values other than $T$.

\end{document}